\documentclass[pre,twocolumn,superscriptaddress]{revtex4}  

\usepackage[paperwidth=210mm,paperheight=297mm,centering,hmargin=2cm,vmargin=2.5cm]{geometry}

\usepackage[pdftex]{graphicx}
\usepackage{amsmath,amsfonts,amssymb}
\usepackage{morefloats}
\usepackage{xcolor} 
\usepackage{tabularx}

\usepackage{newfloat,algcompatible}
\usepackage{etoolbox}
\DeclareFloatingEnvironment[
    fileext=loa,
    listname=List of Algorithms,
    name= Algorithm,
    placement=tbhp,
]{algorithm}

\makeatletter
\renewcommand{\fnum@algorithm}{\small\textbf{\algorithmname~\thealgorithm}}
\makeatother

\definecolor{gris25}{gray}{0.5}
\definecolor{gris5}{gray}{0.7}
\definecolor{gris75}{gray}{0.9}
\definecolor{green2}{rgb}{0,0.8,0}
\definecolor{bluemoi}{rgb}{0.25,0.50 ,0.75} 

\usepackage{hyperref}
\hypersetup{
  bookmarksopen=false,
  pdftitle=" ",
  pdfauthor=" ", 
  pdfsubject=" ", 
  pdftoolbar=true, 
  pdfmenubar=true, 
  pdfhighlight=/O, 
  colorlinks=true, 
  pdfpagemode=UseNone, 
  pdfpagelayout=SinglePage, 
  pdffitwindow=true, 
  linkcolor=bluemoi, 
  citecolor=bluemoi, 
  urlcolor=bluemoi 
}

\makeatletter
\renewcommand{\figurename}{Figure}
\renewcommand{\fnum@figure}{\small\textbf{\figurename~\thefigure}}
\renewcommand{\thefigure}{\arabic{figure}}
\makeatother

\makeatletter
\renewcommand{\tablename}{Table}
\renewcommand{\fnum@table}{\small\textbf{\tablename~\thetable}}
\renewcommand{\thetable}{\arabic{table}}
\makeatother

\begin{document}

\title{A universal model of commuting networks}

\author{Maxime Lenormand}\affiliation{IRSTEA, LISC, 24 avenue des Landais, 63172 AUBIERE, France}
\author{Sylvie Huet}\affiliation{IRSTEA, LISC, 24 avenue des Landais, 63172 AUBIERE, France}
\author{Floriana Gargiulo}\affiliation{IRSTEA, LISC, 24 avenue des Landais, 63172 AUBIERE, France}
\author{Guillaume Deffuant}\affiliation{IRSTEA, LISC, 24 avenue des Landais, 63172 AUBIERE, France}

\begin{abstract} 
We show that a recently proposed model generates accurate commuting networks on 80 case studies from different regions of the world (Europe and United-States) at different scales (e.g. municipalities, counties, regions). The model takes as input the number of commuters coming in and out of each geographic unit and generates the matrix of commuting flows between the units. The single parameter of the model follows a universal law that depends only on the scale of the geographic units. We show that our model significantly outperforms two other approaches proposing a universal commuting model \cite{Balcan2009,Simini2012}, particularly when the geographic units are small (e.g. municipalities).
\end{abstract}

\maketitle

\section*{Introduction}

Billions of people move everyday from home to workplace and generate networks of socio-economic relationships that are the vector of social and economic dynamics such as epidemic outbreaks, information flows, city development and traffic \cite{Ortuzar2011,Balcan2009}. Understanding the essential properties of these networks and reproducing them accurately is therefore a crucial issue for public health institutions, policy makers, urban development, infrastructure planners, etc. \cite{DeMontis2007,DeMontis2010}. This challenge is the subject of an intensive scientific activity (see \cite{Barthelemy2011,Rouwendal2004} for  reviews), in which the analogy of the gravitational attraction inspires a majority of approaches \cite{Wilson1998,Choukroun1975}: the number of commuters between two geographic units (cities, counties, regions...) is supposed proportional to the product of the "masses" of each geographic unit (the population for example) and inversely proportional to a function of the distance between them. Unfortunately, numerous experiments showed that the optimum function and parameter values vary a lot with the case studies \cite{Vries2009,DeMontis2007,DeMontis2010,Fotheringham1981}. This situation is not satisfactory because when one wants to generate a particular commuting network without having the total origin destination matrix of commuting, no practical heuristic is available for choosing the adequate type of function and parameter values. This paper addresses this problem.

We consider a recently proposed model \cite{Gargiulo2012,Lenormand2012}, differentiating itself from the usual gravity law models in two main features:
\begin{itemize}
	\item It takes as input the total number of commuters in and out from each geographic unit. With this starting point, the model focuses directly on the influence of the distance between geographic units on the commuting probability. The model is data demanding, but these data are widely available.
	\item It builds the network progressively, allocating commuters one by one in the different flows, according to probabilities that increase with the number of commuters coming in the destination and decrease with the distance between the origin and destination. These probabilities are updated after each allocation. 
	\end{itemize}

Our model is close to the traditional doubly-constrained gravity model \cite{Wilson1998, Choukroun1975}, but it is more flexible and less data demanding. Indeed, the doubly constrained model and the methods used to solve it require a closed network of commuters: they cannot take into account commuting links outside the considered geographical units. Our individual based stochastic approach overcomes this problem and can deal with the usually available data of total number of commuters in and out of geographic units.

\begin{figure}[b]
	\begin{center}
		\centerline{\includegraphics[width=\linewidth]{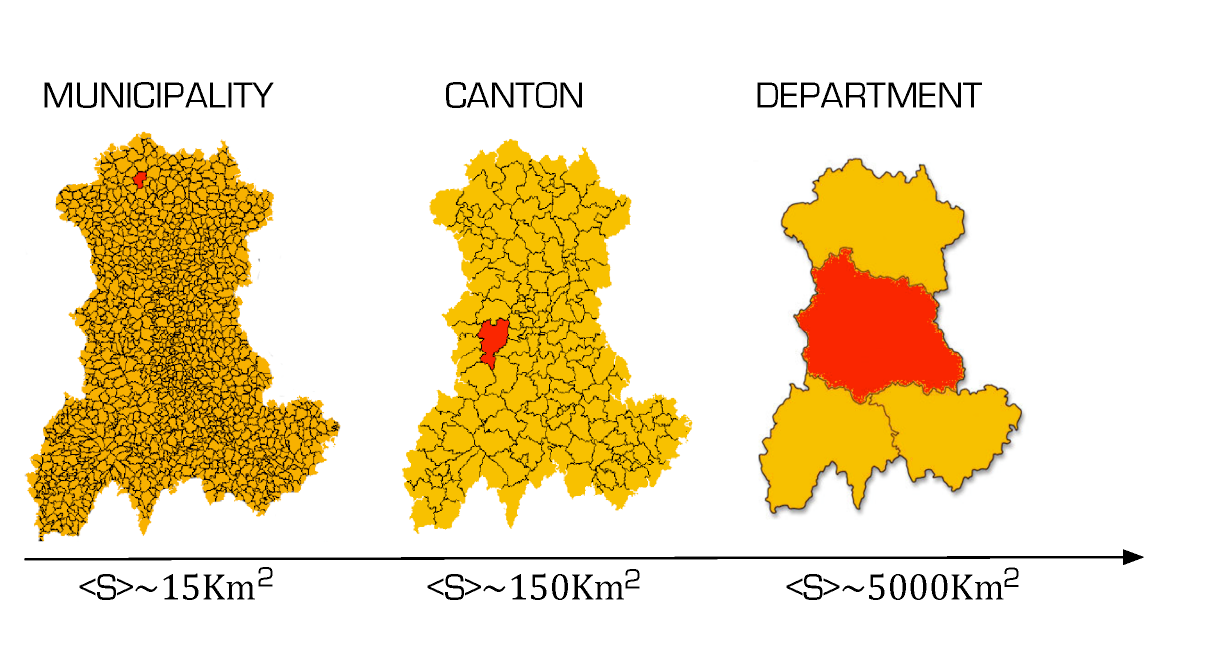}}
		\caption{Three scales of geographic units (Auvergne region, France). \label{Figure1}}
	\end{center}
\end{figure}

We test this model on 80 case-studies with geographic units of different scales. For example in the same case-study  the geographic unit can be either the municipality, the canton or the department, (see an example on Figure \ref{Figure1}). More precisely, the case studies include:  Czech Republic (municipality scale, 1 case-study), France (municipality scale, 34 case-studies), France (canton scale, 15 case-studies including whole France), France (d\'{e}partement scale one case-study (whole France), Italy (municipality scale, 10 case-studies), Italy (province scale, 4 case-studies), USA (county level, 15 case-studies including whole USA). For a detailed description of the datasets see the Appendix.

We show that the single parameter of our model follows a simple universal law that depends only on the average surface of the considered geographic units. This implies that, given the number of commuters in and out of each geographic unit and their average surface, we can derive the whole matrix of flows with a very good confidence. 

Two other approaches \cite{Balcan2009,Simini2012} claim to catch universal properties of commuting networks. We show that our model yields significantly more accurate results, especially for case-studies with small geographic units (e.g. municipalities).

\section*{The model}

We consider the basic double-constrained model setup, without adding any ingredient about the job market characteristics (professions, salary range, etc.). Instead of solving analytically the optimisation problem, we use an individual based procedure that allocates virtual individuals one by one in the different flows between geographic units, according to a probability that is updated after each allocation.

This individual based approach can deal with less constrained data than the doubly-constrained gravity model that requires the total number of commuters in to be equal to the total number of commuters out. In other words the doubly contrained model can only deal with the flows between the considered geographic units; it cannot take into account the commuting links with destinations outside the case study area. This is a problem when only the numbers of commuters in and out the geographic units are available (and not the complete matrix of the commuting flows), because the data do not distinguish between the flows inside and outside the case study area. It is therefore difficult to estimate the correct data to take as input to the doubly-constrained model in this case. Our approach is more flexible and overcomes this difficulty. It does not require that the total number of commuters in and out to be equal (for more details see \cite{Lenormand2012}), hence it can easily use directly the usually available data on the number of commuters in and out of each geographic unit.

Let $s^{out}_i$ and $s^{in}_j$ be respectively the global number of commuters starting from unit $u_i$ and the global number of commuters arriving in unit $u_j$. These numbers are initialised from data and then they are progressively modified by the procedure. More precisely, at each step we select unit $u_i$ such that $s^{out}_i > 0$ at random, and we consider a virtual commuter starting from $u_i$. We draw at random the working place $u_{j^*}$ of this individual among all possible destinations $u_{j}$ according to probabilities $P_{i \rightarrow j}$:
 
\begin{equation}
P_{i \rightarrow j}=\frac{s^{in}_{j} e^{-\beta D_{i j}}}{\sum_{k=1}^{N} s^{in}_k e^{-\beta D_{i k}}}					
\label{probs} 
\end{equation}

where $D_{i j}$ is the Euclidian distance in meter between units $u_i$ and $u_{j}$ (computable from the Lambert or GIS coordinates). Having drawn $u_{j^*}$, we decrement of one $s^{out}_i$ and $s^{in}_{j^*}$. Note that decrementing $s^{in}$ and $s^{out}$ at each step complicates significantly the derivation of an analytical expression of the model. We chose a probability decreasing exponentially with the distance, in accordance with the investigations carried out in \cite{Lenormand2012} and with the literature on commuting network models. The importance of the distance in the commuting choices is embedded in parameter $\beta$: for $\beta\rightarrow 0$ the probability tends to be independent from the distance, while for high values of $\beta$, the probability tends to zero very rapidly when the distance increases, independently from the number of commuters arriving in the units. 

To reduce the border effect (see \cite{Lenormand2012}), we consider the job-search basin in an extended (EXT) area, composed by the $n$ residential units and $m$ units surrounding the area. Thus, we have $n$ units which are commuting origins and $N=n+m$ units that are commuting destinations. The generated network is saved in matrix $\tilde T\in M_{n\times N}(\mathbb{N})$ where each entry $\tilde T_{i j}$ represents the number of commuters between units $u_i$ and $u_j $. The algorithm is summarized in Algorithm \ref{Algorithm1}.

\begin{algorithm}[!ht]
    {\hrulefill}
	\vspace*{-0.3cm}
	\caption{Commuting generation model}
		{\vspace*{-0.25cm}\hrulefill}
	\label{Algorithm1}
			\begin{algorithmic}
	  	\REQUIRE $D \in \mathrm{M}_{n\times N}(\mathbb{R})$, $s^{in} \in \mathbb{N}^{N}$, $s^{out} \in \mathbb{N}^n$, $\beta \in \mathbb{R}_+$
			\ENSURE $\tilde T \in \mathrm{M}_{n\times N}(\mathbb{N})$
			\STATE $\tilde T_{i j} \leftarrow 0$
			\WHILE {$\sum_{k =1}^n s^{out}_k > 0$}
				\STATE Pick at random $i  \in |[1,n]|$, such that $s^{out}_{i} \neq 0$
				\STATE Pick at random $j$ from  $|[1,N]|$ 
				\STATE \ \ \ \ with a probability $P_{i \rightarrow j}$
				\STATE $\tilde T_{i j} \leftarrow \tilde T_{i j}+1$		
				\STATE $s^{in}_j \leftarrow s^{in}_j-1$ 
				\STATE $s^{out}_i \leftarrow s^{out}_i-1$	
			\ENDWHILE
		\end{algorithmic}
	{\vspace*{-0.2cm}\hrulefill}		
\end{algorithm} 

\section*{Results}

\subsection*{A universal law ruling parameter $\beta$}

The model depends on a single parameter ruling the importance of the distance in commuting choice. We show that this parameter can be derived as a function of the scale of the problem, independently from the socio-geographical location of the case study area. This opens the possibility to reconstruct the commuting flows (origin-destination matrix) when they are not provided.

\begin{figure}
	\begin{center}
		\centerline{\includegraphics[width=\linewidth]{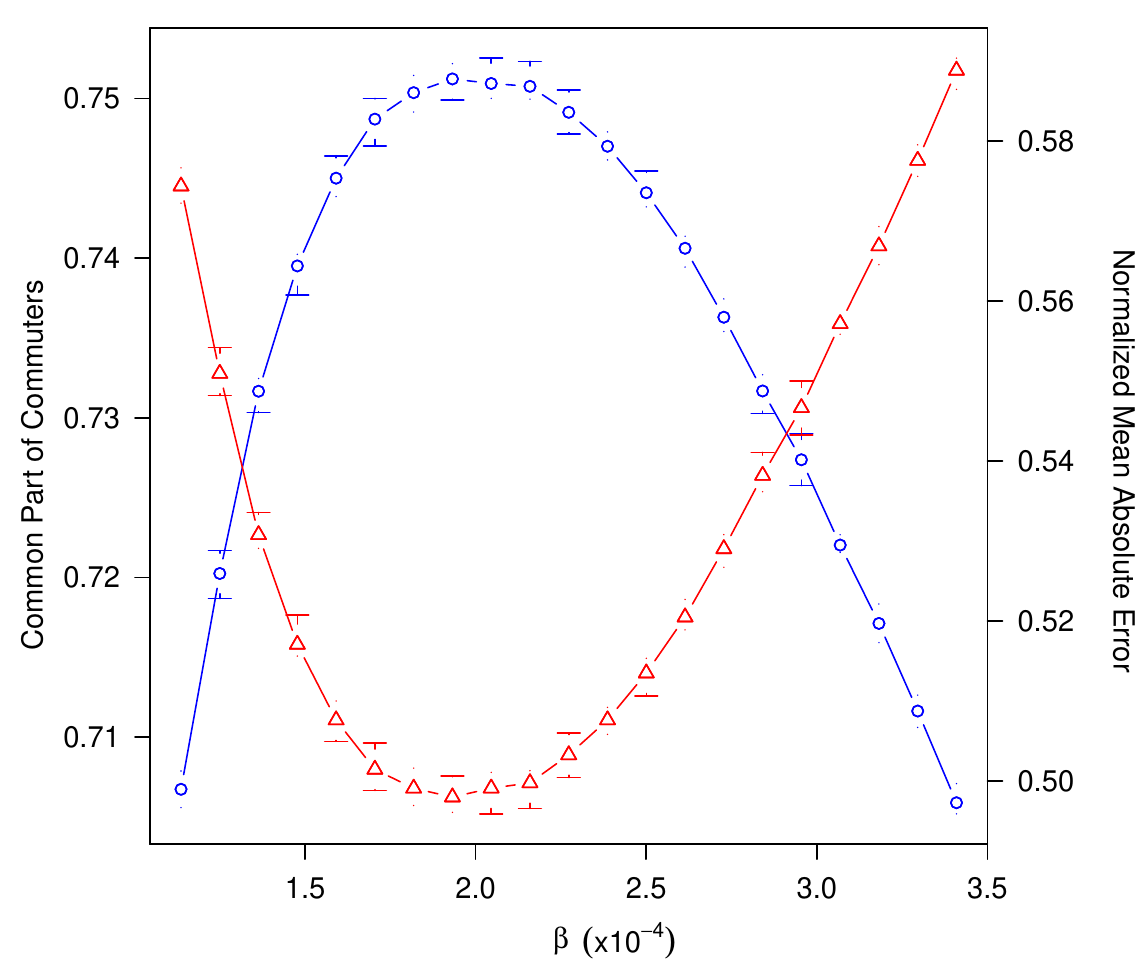}}
		\caption{Plot of the average CPC (blue circle) and the average NMAE (red triangle) in term of $\beta$ for 10 replications of the model for the Auvergne case study (FR1). The error bars represent the minimum value and maximum value obtain over the 10 replications. \label{Figure2}}
	\end{center}
\end{figure}

We calibrated parameter $\beta$ by maximising the common part of commuters (CPC), based on the S{\o}rensen index \cite{Sorensen1948}. 

\begin{equation}
	CPC(T,\tilde {T})=\frac{2 NCC(T,\tilde {T})}{NC(T) + NC(\tilde {T})}
\end{equation}
with:
\begin{equation}
	NCC(T,\tilde {T})=\sum_{i=1}^{n}\sum_{j=1}^{n}\min(T_{ij}\tilde{T}_{ij})
\end{equation}

and

\begin{equation}
	NC(T)=\sum_{i=1}^{n}\sum_{j=1}^{n} T_{ij}
\end{equation}
 
where $T$ is the observed origin-destination matrix and $\tilde{T}$ is the simulated one. This is a similarity measure based on the S{\o}rensen index in ecology computing which part of the commuting flows is correctly reproduced, on average, by the simulated network. It varies between 0, when no agreement is found, and 1, when the two networks are identical. We priviledged this indicator because of its direct interpretation. Indeed, when $NC(T) \simeq NC(\tilde {T})$ (it is the case for our model), the CPC represents the percentage of commuting connection correctly located (i.e. with the right pair origin - destination). Moreover, we tested on all case studies that the results obtained with the MAE, the RMSE or CPC are equivalent (see the Appendix for more details). We have also shown in \cite{Gargiulo2012, Lenormand2012} that the value of $\beta$ yielding the maximum CPC also yields the maximum similarity between observed and simulated commuting distance distributions. As an example on the $FR1$ case study, Figure \ref{Figure2} shows that the same $\beta$ value maximizes the CPC and minimizes the MAE. In this Figure we can also note that the CPC is very sensitive to $\beta$ and that its value does not vary much with the different replicas of the stochastic solving process.

Moreover, in order to have an idea of the improvement of the model compared with complete randomness, we have computed the CPC of a random model where the probabilities presented in Equation (\ref{probs}) are uniform ($P_{i \rightarrow j}= \frac{1}{n}$, where $n$ is the number of units). As shown on the Figure \ref{Figure3} we obtained an average CPC around 0.1. For our model, the CPC is always higher than 0.7 with an average around 0.8, which can be interpreted as 70 to 80 \% of correctly predicted commuting connections.

\begin{figure}[b]
	\begin{center}
		\includegraphics[width=\linewidth]{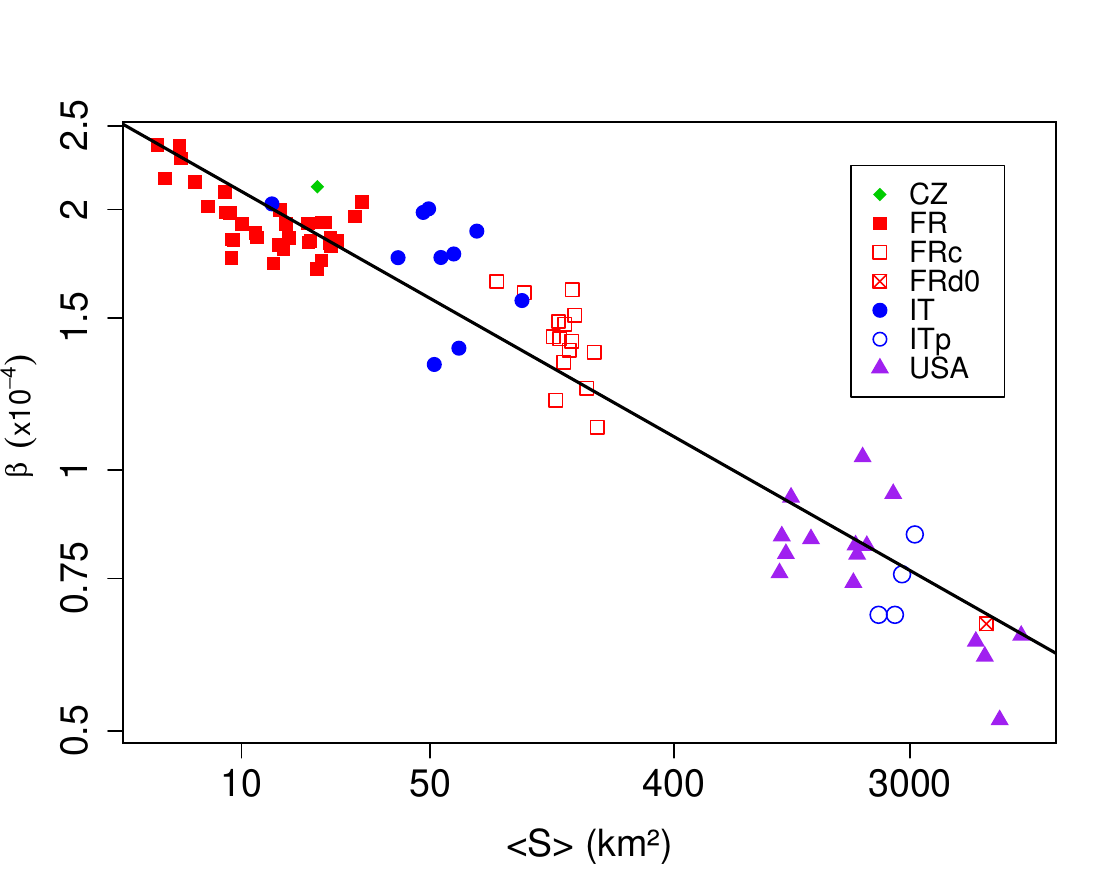}
		\caption{Log-log scatter plot of the calibrated $\beta$ values in terms of average surface of the geographic units for 80 case-studies; the line represents the regression line predicting $\beta$. \label{Figure4}}
	\end{center}
\end{figure}

\begin{figure*}
	\begin{center}
		\centerline{\includegraphics[scale=0.8]{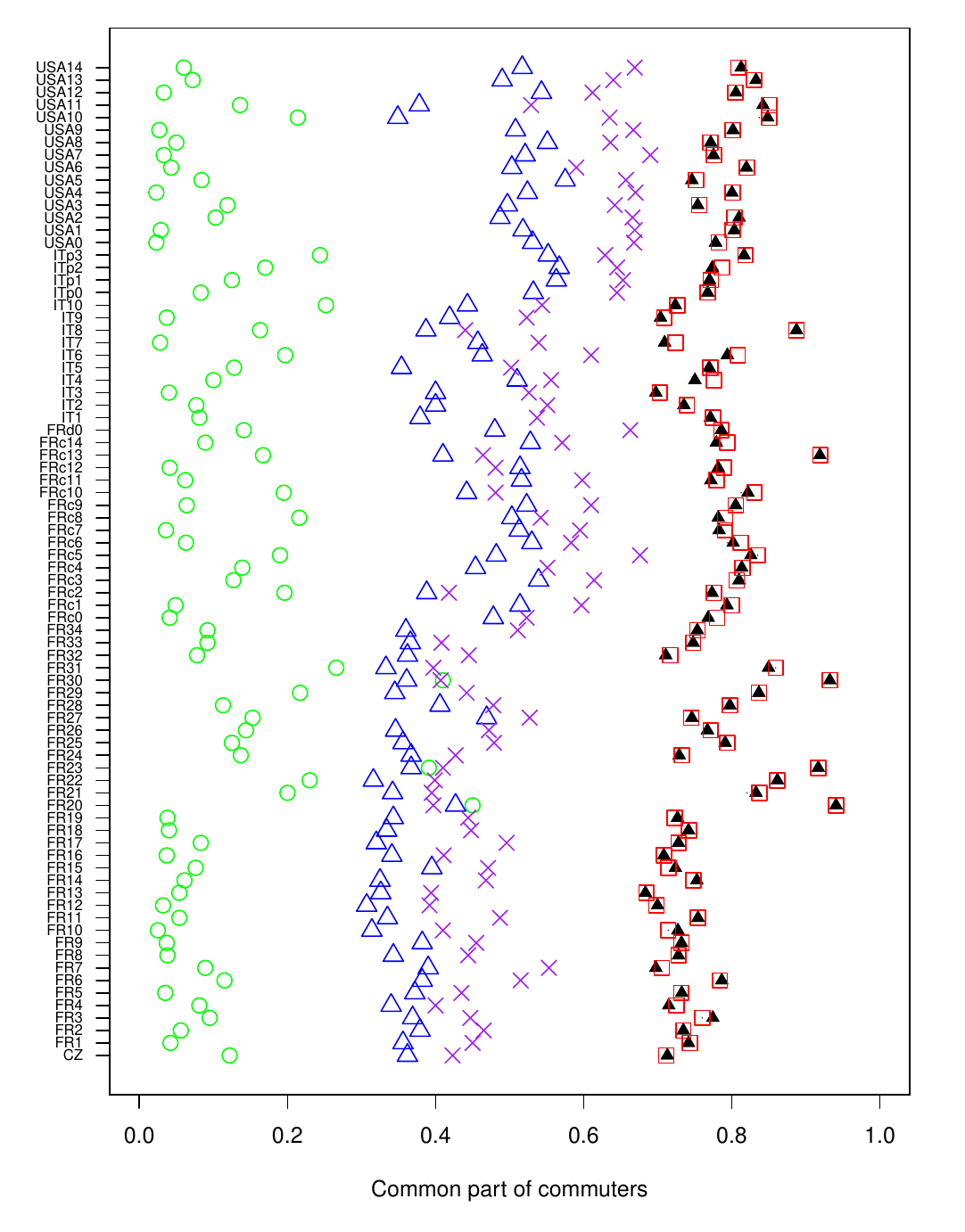}}
		\caption[Common part of commuters (CPC) for the 80 case-studies]{Common part of commuters (CPC) for the 80 case-studies. The red squares represent 
		the CPC obtained with the value of $\beta$ optimised from data on the case-study network. Black plain triangles represent the average CPC obtained 
		with $\beta$ values estimated with the rule linking $\beta$ and the average surface of the units obtain with the cross-validation; Dark bars 
		represent the minimum and the maximum CPC obtained with the estimated $\beta$ but in most cases they are too close to the average to be seen. The 
		green circles represent the CPC obtained with the random model. The blue triangles represent the CPC obtained with the radiation model. The purple 
		crosses represent the CPC obtained with the modified version of the radiation model.
		\label{Figure3}}
	\end{center}
\end{figure*}

Our goal is to derive the value of $\beta$ from some easily available global characteristics of the case-study, giving the possibility to reconstruct the commuting flows when they are not available. Figure \ref{Figure4} gives strong evidence of such a universal relation. 

The x-axis represents the average surface of the geographic units of the case-study ($\langle S\rangle$ in logarithm scale) and the y-axis the optimal $\beta$ value (in logarithm scale). The linear regression in the log-log plane shows a simple relation:

\begin{equation}
\beta=\alpha\langle S\rangle^{-\nu}\label{loglogcalib}
\end{equation}

with $\alpha=3.15\cdot 10^{-4}$ and $\nu=0.177$. $\alpha$ corresponds to the $\beta$ value for the unitary surface 1 $km^2$. The high value of the adjusted $R^2=0.92$ confirms the quality of the linear model. We observe that $\beta$ decreases with the average surface of the units $\langle S\rangle$, meaning that, when $\langle S\rangle$ is small (e.g. for municipalities in France) the distance is more important in the commuting choice than when $\langle S\rangle$ is large (e.g. for regions or counties).

We now evaluate the robustness of our estimation of $\alpha$ and $\nu$ using a common statistical procedure: the cross-validation. The cross-validation aims at evaluating the potential error of using the $\beta$ value derived from the regression model intsead of deriving this value by optimisation for a new case study. This procedure repeats a large number of times the following steps: define a sub-sample of the total sample of case studies, derive a regression model of $\beta$ from this sub-sample, for each case study that do not belong to the sub-sample, derive $\beta$ from this regression model and compare the corresponding CPC with the value of $\beta$ directly calibrated on the complete origin - destination  data. The dataset (including 80 case-studies) is randomly cut into two sets, called the training set (comprising 53 case-studies) and the test set (composed of 27 case-studies). We build a regression model on the training set, providing $\alpha$ and $\nu$, from which we derive estimates of $\beta$ for each of the 27 case-studies of the testing set. We have 27 estimations of $\beta$ using the relation \ref{loglogcalib} where $\alpha$ and $\nu$ are obtained from the random sub-sample of 53 case-studies. We repeat this process $10,000$ times obtaining 270,000 estimations of $\beta$ (uniformly distributed over the 80 case-studies) corresponding to about $\frac{270,000}{80}=3,375$ estimations of $\beta$ for each case study. Then we calculate the average, minimum and maximum CPC for each of these values of $\beta$, and we compare them with the CPC obtained with value of $\beta$ directly calibrated on the data. 

Figure \ref{Figure3} shows, for each case-study, the CPC associated with the calibrated $\beta$, the average CPC obtained with the $\beta$ values estimated from the cross-validation and the confidence interval defined by the minimum and the maximum values (but it is too small to be seen in most cases). The CPC obtained with the calibrated $\beta$ value (black triangle) is almost the same as the average CPC obtained with the estimated $\beta$ in most cases (red square). Globally, we can conclude that the $\beta$ estimated with the log-linear model and the calibrated $\beta$ lead to very similar CPCs and also very similar MAE and the RMSE as shown in the Appendix. The method appears therefore fairly robust and this gives confidence for using it with the value of $\beta$ derived from our loglog regression in new cases studies. 

\begin{figure*}
				\includegraphics[width=\linewidth]{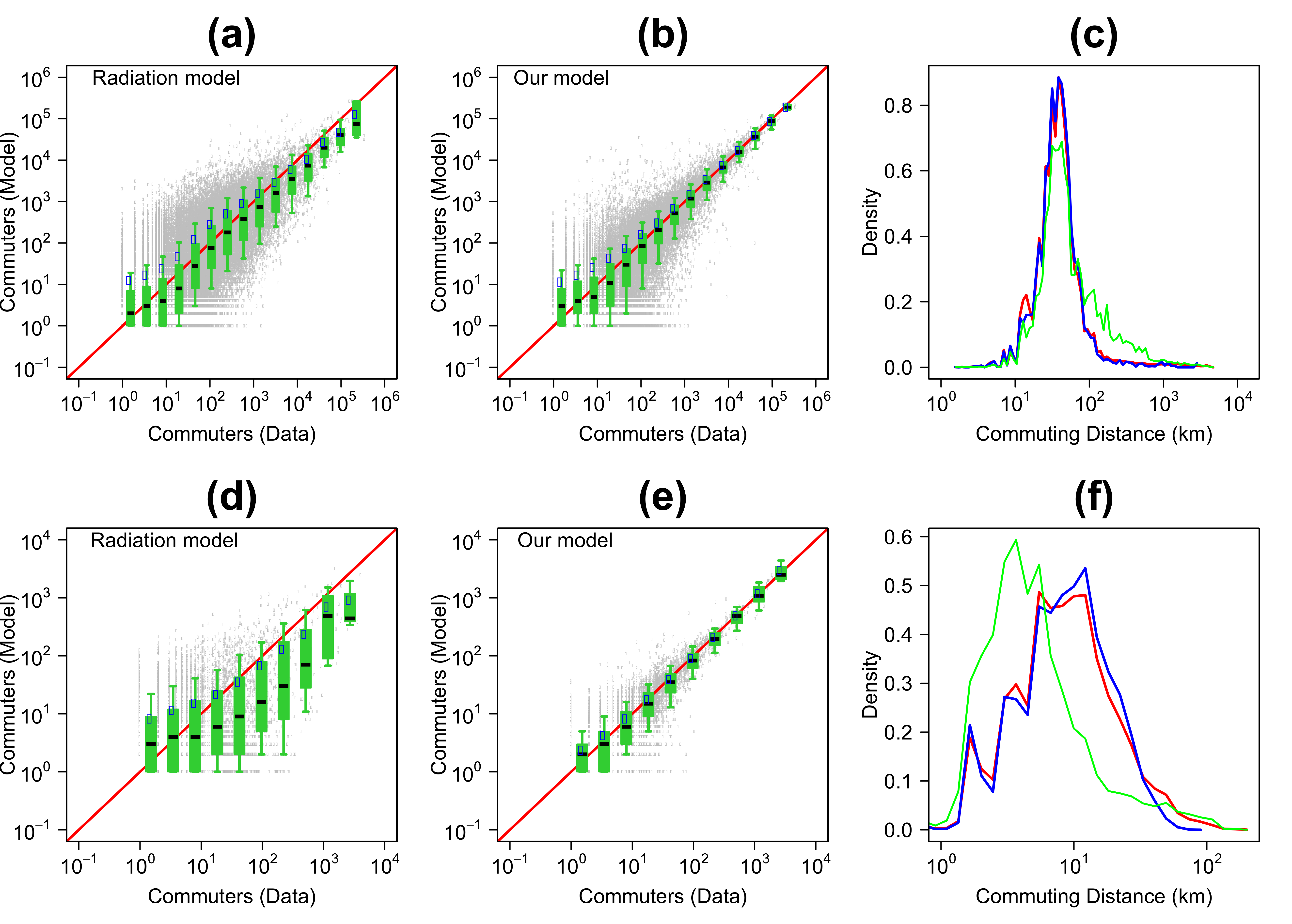}		
				\caption[Comparing the predictions of the radiation model with ours]{Comparing the predictions of the radiation model with ours for two case 
				studies, the first row ((a)-(c)) for USA0 (USA at 
				county scale) and the second row ((d)-(f)) for FR1  (Auvergne region, France at municipality scale). Plots (a), (b) , (d) and (e): Comparison 
				between the observed (Census) and the simulated (model) non-zero flows. Grey points are the scatter plot for each pair of units. 
				The boxplots (D1, Q1, Q2, Q3 and D9) represent the distribution of the number of simulated travelers in different bins of number of observed 
				travelers. The blue circles represent the average number of simulated travelers in the different bins. Plots (c) and (f): Commuting distance 
				distributions (km) (i.e. Probability for a commuters of the region to commut at a distance $d$). The blue line represents the observed data, 
				the red one the results of our model and the green one the results of the radiation model. \label{Figure5}} 
\end{figure*}

\subsection*{Comparaison with other universal derivations of commuting networks}

Two other different approaches, \cite{Balcan2009} and \cite{Simini2012}, claim also to provide a universal derivation of commuting networks. 
The objective of \cite{Balcan2009} is to generate a worldwide commuting network, and the model must deal with the wide variety of populations and surfaces of geographic units for which the data are available. To solve this difficulty, the authors project these data on ad-hoc units defined with a Voronoi diagram. They define their basic unit as a cell approximately equivalent to a rectangle of 25 x 25 kilometers along the Equator. This allows them to calibrate their model because a unit is the same object whatever the country. This is an interesting solution for generating a world-wide commuting network but it leads to an average commuting distance of 250 km which is much larger than the average distance of daily commuting. For example for the USA case study the average distance of daily commuting is about 68 km for the observed network and about 64 km for the simulated network obtained with our algorithm. For the Auvergne (France) case study at municipality scale the average distance of daily commuting is about 12 km for the observed network and about 11 km for the simulated one. 

In the radiation model, proposed in \cite{Simini2012}, the commuting flow between two geographic units is a function of the cumulated population in a circle at the distance between the two units. The model has an elegant analytical solution and the average flow $T_{ij}$ from unit $u_i$ to unit $u_j$ can be approximated by 
\begin{equation}
\langle T_{ij}\rangle = \left(m_{i}\frac{P_{c}}{P}\right)\frac{m_{i}n_{j}}{\left(m_{i}+s_{ij}\right)\left(m_{i}+n_{j}+s_{ij}\right)}
 \end{equation}
where $m_{i}$ and $n_{j}$ are respectively the population of units $u_i$ and $u_j$, $P_c$ is the total number of commuters and $P$ is the total population in the case-study region, and $s_{ij}$ the total population in the circle of radius $r_{ij}$ centred at $u_i$ (excluding the source and destination population). 

We implemented their analytical approximation and reproduced the graphs presented in their paper. Figure \ref{Figure5} shows the comparison between the radiation model and ours in the US for inter-county commuting and in the French Auvergne region for inter-municipality commuting (see the Appendix for more examples). We observe that in both cases our approach yields significantly better results.
Moreover, as shown on Figure \ref{Figure3}, the average CPC for the radiation model on all the case studies is around 0.4, and lower for all case studies than the one obtained with our approach. 

However, it should be reminded that our model uses more specific data (total number of commuters in and out of each geographic unit) than the radiation model, hence one could expect our results to be more accurate. Therefore, to be fair with the radiation model we implemented a modified version of this model using the number of out and in commuters of each units. This new approximation is presented in Equation \ref{Sim2} where $s_{ij}$ the total number of in-commuters in the circle of radius $r_{ij}$ centred at $u_i$ (excluding the source and destination).

\begin{equation}
	\langle T_{ij}\rangle = s^{out}_i \frac{s^{out}_i s^{in}_j}{\left(s^{out}_i + s_{ij}\right)\left(s^{out}_i + s^{in}_j + s_{ij}\right)}\label{Sim2}
\end{equation}

As shown on Figure \ref{Figure3}, this new model reaches an average CPC around 0.5 which is higher than the original radiation model but still significantly lower than the results obtained with our model. Using the MAE and the RMSE leads to the same conclusions (see the Appendix for more details).

\section*{Discussion}

The power law of our model's single parameter $\beta$ with the average area of the case study geographic units, is surprising to us because of the high variety in our case studies in terms of scale, number of units, number of commuters and surface areas. For instance the Auvergne region in France is rural with a population density of about 50 hab./km$^2$ whereas the New York City region is very urban with a population density of about 6500 hab./km$^2$. As far as we know, this is the first time that a single model is shown to fit such diverse group of datasets.

We show that our approach outperforms the radiation model and that the difference of input data plays a minor role in this superiority. This superiority is not due to our particular treatment of the border effects either. Indeed, we could check our approach outperforms the radiation model also on particular case studies (e.g. on islands such as Corsica) where this border effect does not play. We can conclude that the accuracy of our model comes from a proper use of the number of commuters in and out of each geographic unit and an adequate choice of the function of the distance. 

The results of the cross validation procedure give a good confidence in the robustness of this law. However, we have to admit that, despite their diversity, our 80 case studies come all from western industrialised countries. Therefore it will be important to check the validity of our law on case studies coming from other continents and less industrialised countries. Moreover, we use a very rough approximation of the distance between the geographic units with the Euclidian distance between the unit centroids. More accurate approximations of this distance would certainly improve the results. Finally, we also intend to apply our approach to commuting networks inside urban areas because many cities of the world show an impressive growth and an increasing part of commuting takes place within them \cite{Roth2011}. An important issue in our perspective is to check if our law holds at this scale.

\section*{Acknowledgments}

This publication has been funded by the Prototypical policy impacts on multifunctional activities in rural municipalities collaborative project, European Union 7th Framework Programme (ENV 2007-1), contract no. 212345. The work of the first author has been funded by the Auvergne region.

\bibliographystyle{unsrt}  
\bibliography{Network}

\newpage

\makeatletter
\renewcommand{\fnum@table}{\small\textbf{\tablename~S\thetable}}
\makeatother
\setcounter{table}{0}

\makeatletter
\renewcommand{\fnum@figure}{\small\textbf{\figurename~S\thefigure}}
\makeatother
\setcounter{figure}{0}

\section*{APPENDIX}

\subsection*{Data presentation}

\subsubsection*{Datasets}

Commuting data are usually provided by statistical offices in the form of origin-destination tables. We analyzed $80$ case studies from $7$ differents datasets and $4$ different country (described in Table S\ref{Datasets}). In these appendices we called outside (Out.) the $m$ units surrounding the area. 

\subsubsection*{Distance}

The distances between units are Euclidean, computed using the Lambert coordinates or the latitude/longitude of the centroid of the units. 

\subsubsection*{Case studies}

We define two types of case studies: from administrative regions and from aggregation of small administrative units around a randomly chosen point. Each case study is composed of a region and an outside (the units surrounding the region at a reasonable distance).
 
\begin{figure*}
  \includegraphics[width=\linewidth]{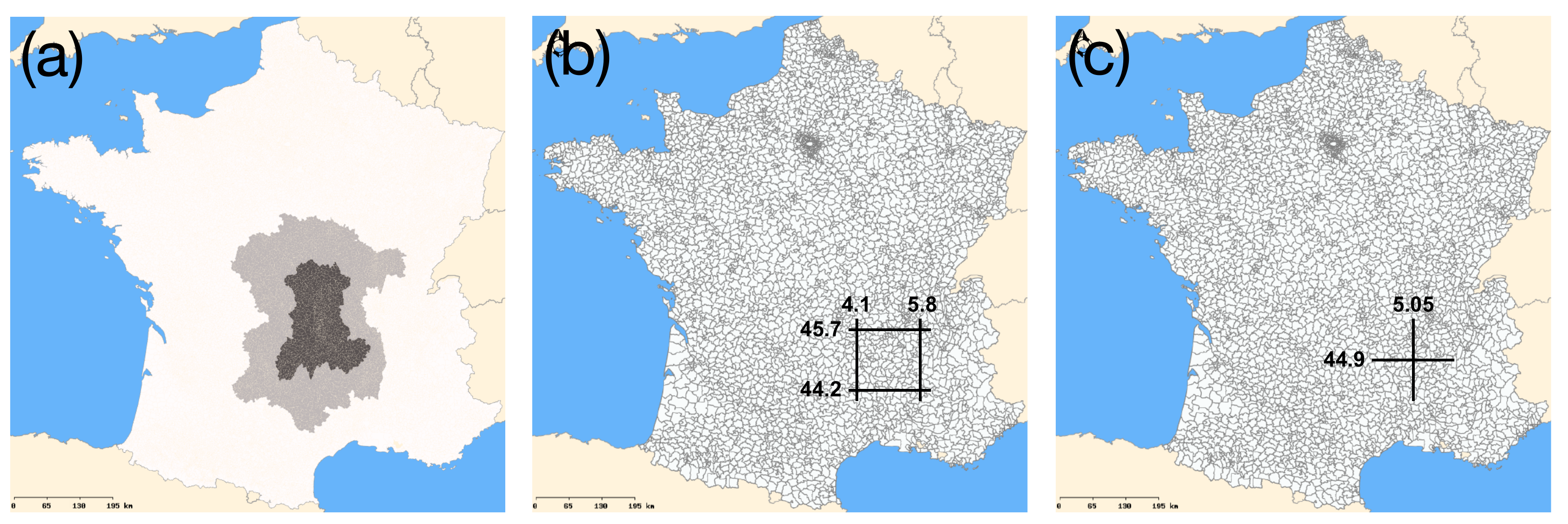}
  \caption{Maps to illustrate the build process regions. (a) Administrative; (b) starting point of aggregation and (c) limits of aggregated units. \textit{\scriptsize{Base maps source: Cemagref - DTM - D\'{e}veloppement Informatique Syst\`{e}me d'Information et Base de Donn\'{e}es : F.Bray \& A.TorreIGN (G\'{e}ofla\textsuperscript{{\fontsize{5}{5}\textregistered}},2007).}}}
  \label{FigureS1}
\end{figure*}

To build a case study from an admistrative region, we select an administrative region (For example the Auvergne region represented by the dark grey region in Figure \ref{FigureS1}a) and to build the outside we select all the units surrounding the region at a reasonable distance (For the Auvergne region example, the outside is represented by the light grey region in Figure \ref{FigureS1}a). 

To build a case study by aggregation of units, firstly, we define the number of desired units and we draw at random a latitude and a longitude (For example the point represented in Figure \ref{FigureS1}b). In a second time we gradually increase the area of a square with as center the starting point until the desired number of units is obtained (Figure \ref{FigureS1c}). To build the outside we select all the units surrounding the defined set of units at a reasonable distance or all the remaining units in the country (it depends of the number of units). 

The case studies with an identifier with a $0$, for example $FRc0$, are complete network of the country without outside. Indeed, we have no data for the surrounding countries. When we consider a region in the country we can determine the outside as the units surrounding the region. When we consider as a region the whole country we can't determine an outside, it is the case for $FRc0$ (all the Cantons of France), $Frd0$ (all the D\'epartements of France), $Itp0$ (all the Provincias of Italy) and $USA0$ (all the counties of USA).

\subsubsection*{Source}

The $3$ French datasets are measured for the 1999 French Census by the French Statistical Institute, $INSEE$. They were kindly made available by the Maurice Halbwachs Center. The $2$ Italian datasets are measured for the 2001 Italian Census by the National Institute for Statistics, $ISTAT$. 

\begin{table*}
\caption{Description of the case studies}
\label{DescriptReg}
\begin{center}
\fontsize{7}{7}\selectfont
\begin{tabular}{>{\centering}m{1cm}>{\centering}m{1.5cm}>{\centering}m{1.5cm}>{\centering}m{1cm}>{\centering}m{1cm}>{\centering}m{1.5cm}>{\centering}m{2cm} m{2cm}<{\centering}}
\hline
\textbf{Case study} & \textbf{Number of units (area)} &  \textbf{Number of units (outside)} & \textbf{Surface (km$^2$)} & \textbf{Average unit surface (km$^2$)} &\textbf{Standard deviation unit surface (km$^2$)} & \textbf{Observed number of commuters (Area)} & \textbf{Estimated number of commuters (Area)}\\
\hline
 & & & & & & & \\
CZ	&	43	&	630	&	35369	&	822.54	&	703.23	&	6585	&	6847	\\
FR1	&	1310	&	3463	&	26013	&	19.86	&	12.49	&	261822	&	262452	\\
FR2	&	1269	&	1447	&	27208	&	21.44	&	16.14	&	608587	&	613363	\\
FR3	&	419	&	2809	&	5762	&	13.75	&	8.46	&	90456	&	76829	\\
FR4	&	903	&	3081	&	8280	&	9.17	&	9.55	&	409661	&	402565	\\
FR5	&	2296	&	2835	&	41309	&	17.99	&	21.30	&	679639	&	657095	\\
FR6	&	261	&	3124	&	5175	&	19.83	&	10.46	&	52921	&	48681	\\
FR7	&	185	&	1859	&	5167	&	27.93	&	18.71	&	9474	&	8981	\\
FR8	&	1464	&	2467	&	25810	&	17.63	&	12.94	&	333045	&	333540	\\
FR9	&	1842	&	4718	&	39151	&	21.25	&	14.76	&	514461	&	529535	\\
FR10	&	3020	&	3845	&	45348	&	15.02	&	15.74	&	502326	&	494946	\\
FR11	&	747	&	3169	&	16942	&	22.68	&	14.15	&	118508	&	117217	\\
FR12	&	1786	&	3317	&	16202	&	9.07	&	7.46	&	239931	&	236314	\\
FR13	&	1420	&	3536	&	12317	&	8.67	&	5.64	&	396800	&	402128	\\
FR14	&	433	&	3914	&	6211	&	14.34	&	12.41	&	30175	&	28729	\\
FR15	&	515	&	3808	&	5874	&	11.41	&	9.54	&	76519	&	72896	\\
FR16	&	2339	&	3067	&	23547	&	10.07	&	7.51	&	505807	&	507812	\\
FR17	&	260	&	1814	&	5565	&	21.40	&	13.15	&	17310	&	17071	\\
FR18	&	1545	&	3046	&	27367	&	17.71	&	15.78	&	354824	&	354566	\\
FR19	&	1948	&	1983	&	25606	&	13.14	&	12.94	&	333045	&	329908	\\
FR20	&	36	&	1245	&	176	&	4.89	&	3.28	&	193236	&	182808	\\
FR21	&	262	&	1543	&	2284	&	8.72	&	6.62	&	226205	&	206624	\\
FR22	&	185	&	1707	&	1246	&	6.74	&	3.83	&	143938	&	124185	\\
FR23	&	47	&	1234	&	245	&	5.21	&	3.03	&	143586	&	121474	\\
FR24	&	377	&	2283	&	3525	&	9.35	&	7.44	&	160294	&	157123	\\
FR25	&	195	&	2338	&	3718	&	19.07	&	17.66	&	26576	&	24975	\\
FR26	&	547	&	449	&	4116	&	7.52	&	15.87	&	59709	&	61324	\\
FR27	&	163	&	353	&	4299	&	26.37	&	27.53	&	145995	&	148922	\\
FR28	&	327	&	2788	&	4781	&	14.62	&	9.76	&	134048	&	130910	\\
FR29	&	102	&	2031	&	609	&	5.97	&	4.21	&	22520	&	20549	\\
FR30	&	40	&	783	&	236	&	5.90	&	4.28	&	139181	&	125542	\\
FR31	&	196	&	1597	&	1804	&	9.20	&	6.04	&	188855	&	165505	\\
FR32	&	463	&	2588	&	5229	&	11.29	&	8.03	&	50505	&	51413	\\
FR33	&	433	&	2728	&	6004	&	13.87	&	9.07	&	69377	&	63078	\\
FR34	&	286	&	2088	&	5857	&	20.48	&	13.36	&	38141	&	37197	\\
FRc0	&	3646	&	0	&	540241	&	171.72	&	99.90	&	12193161	&	12193161	\\
FRc1	&	1062	&	2584	&	173797	&	163.65	&	91.23	&	2229003	&	2265247	\\
FRc2	&	523	&	3123	&	58366	&	111.60	&	114.44	&	3892543	&	3922481	\\
FRc3	&	226	&	3420	&	33041	&	146.20	&	70.56	&	548048	&	558086	\\
FRc4	&	160	&	3486	&	25044	&	156.52	&	75.47	&	320432	&	323169	\\
FRc5	&	55	&	3591	&	7847	&	142.67	&	71.64	&	61761	&	60285	\\
FRc6	&	869	&	2777	&	131174	&	150.95	&	96.62	&	1995302	&	1983097	\\
FRc7	&	2088	&	1558	&	351073	&	168.14	&	94.18	&	4459338	&	4523902	\\
FRc8	&	100	&	3546	&	20246	&	202.46	&	161.41	&	307744	&	316592	\\
FRc9	&	600	&	3046	&	113905	&	189.84	&	103.57	&	1078183	&	1095993	\\
FRc10	&	302	&	3344	&	26627	&	88.17	&	77.64	&	1306425	&	1274670	\\
FRc11	&	906	&	2740	&	142619	&	157.42	&	100.21	&	2324444	&	2358580	\\
FRc12	&	1500	&	2146	&	250676	&	167.12	&	99.00	&	3224586	&	3284517	\\
FRc13	&	32	&	3614	&	6653	&	207.91	&	145.33	&	11959	&	10634	\\
FRc14	&	506	&	3140	&	75603	&	149.41	&	85.63	&	1311912	&	1331984	\\
FRd0	&	94	&	0	&	540250	&	5747.35	&	1957.11	&	3548178	&	3548178	\\
IT1	&	377	&	0	&	24090	&	63.90	&	61.89	&	225351	&	225351	\\
IT2	&	395	&	201	&	24157	&	61.16	&	77.51	&	409889	&	408692	\\
IT3	&	1002	&	2020	&	54918	&	54.81	&	71.37	&	1235378	&	1193338	\\
IT4	&	201	&	507	&	14964	&	74.45	&	82.42	&	246609	&	248562	\\
IT5	&	204	&	1005	&	10567	&	51.80	&	55.68	&	279014	&	272310	\\
IT6	&	51	&	506	&	5582	&	109.45	&	101.52	&	57446	&	51211	\\
IT7	&	2000	&	4001	&	98693	&	49.35	&	60.97	&	2849914	&	2812238	\\
IT8	&	186	&	1023	&	2412	&	12.97	&	15.25	&	316602	&	286285	\\
IT9	&	1510	&	4004	&	71167	&	47.13	&	58.08	&	1703944	&	1702002	\\
IT10	&	705	&	3008	&	26809	&	38.03	&	41.62	&	401998	&	403307	\\
ITp0	&	99	&	0	&	277220	&	2800.20	&	1619.86	&	1567576	&	1567576	\\
ITp1	&	50	&	49	&	131773	&	2635.45	&	1401.23	&	742229	&	727038	\\
ITp2	&	30	&	69	&	93666	&	3122.21	&	1599.56	&	266696	&	272316	\\
ITp3	&	20	&	79	&	45854	&	2292.72	&	1128.38	&	264824	&	259988	\\
USA0	&	3108	&	0	&	8070785	&	2596.78	&	3437.29	&	34077841	&	34077841	\\
USA1	&	1015	&	2093	&	1876151	&	1848.42	&	916.86	&	5855813	&	5902784	\\
USA2	&	103	&	3005	&	101411	&	984.57	&	341.47	&	527136	&	535608	\\
USA3	&	54	&	3054	&	306284	&	5671.93	&	4488.99	&	604043	&	597371	\\
USA4	&	2011	&	1097	&	4169235	&	2073.21	&	1786.40	&	14767588	&	14926726	\\
USA5	&	202	&	2906	&	404093	&	2000.46	&	1994.32	&	8789633	&	8893748	\\
USA6	&	504	&	2604	&	949238	&	1883.41	&	1041.57	&	2125887	&	2155981	\\
USA7	&	806	&	2302	&	4234740	&	5254.02	&	5626.18	&	5003104	&	5099317	\\
USA8	&	352	&	2756	&	2723212	&	7736.40	&	7741.02	&	4147054	&	4234376	\\
USA9	&	1507	&	1601	&	2877429	&	1909.38	&	1517.28	&	10099598	&	10234438	\\
USA10	&	13	&	3095	&	14123	&	1086.37	&	343.73	&	58212	&	53513	\\
USA11	&	32	&	3076	&	205989	&	6437.17	&	4105.95	&	22496	&	24085	\\
USA12	&	1004	&	2104	&	1292835	&	1287.68	&	563.79	&	9704950	&	9735646	\\
USA13	&	207	&	2901	&	207785	&	1003.79	&	352.24	&	1307774	&	1326018	\\
USA14	&	301	&	2807	&	312955	&	1039.72	&	394.71	&	2054878	&	2085408	\\
\hline   
\end{tabular}
\end{center}
\end{table*}

\subsection*{Results with standard indicators of error}

We computed the results with standard indicators of error.
\begin{itemize}
	\item The Normalized Mean Absolute Error:
\begin{equation}
	NMAE(T,\tilde {T})=\frac{\sum_{i=1}^{n}\sum_{j=1}^{n}|T_{ij}-\tilde{T}_{ij}|}{\sum_{i=1}^{n}\sum_{j=1}^{n} T_{ij}}
\end{equation}
 \item Normalized Root Mean Square Error:
 \begin{equation}
	NRMSE(T,\tilde {T})=\frac{\sqrt{\sum_{i=1}^{n}\sum_{j=1}^{n}(T_{ij}-\tilde{T}_{ij})^2}}{\sum_{i=1}^{n}\sum_{j=1}^{n} T_{ij}}
\end{equation}
\end{itemize}

\begin{table*}
	\caption{Presentation of the datasets. * Data are available online at \url{http://www.czso.cz/eng/redakce.nsf/i/home} and ** Data are available online at 
	\url{http://www.census.gov/geo/www/gazetteer/places2k.html}}
	\label{Datasets}
			\begin{tabular}{|>{\centering}m{1.5cm}|>{\centering}m{1.5cm}|>{\centering}m{1.5cm}|>{\centering}m{1.7cm}|>{\centering}m{2.2cm}|>{\centering}m{2cm}|>{\centering}m{1cm}| 	
			m{1cm}<{\centering}|}
			  \hline
				\textbf{Dataset} & \textbf{Country} & \textbf{Case Study} & \textbf{Distance} & \textbf{Region} & \textbf{Scale} & \textbf{Year} & \textbf{Source}\\
				\hline
				1 & Czech Republic & CZ & Latitude Longitude & Administrative & Municipality & 2001 & * \\
				\hline
				2 & France & FR1 - FR34 & Lambert & Administrative & Municipality & 1999 & INSEE\\
				\hline
				3 & France & FRc0 - FR14 & Latitude Longitude &  Arbitrary aggregation & Canton & 1999 & INSEE\\
				\hline
				4 & France & FRd0 & Latitude Longitude & Administrative & D\'epartement & 1999 & INSEE\\
				\hline
				5 & Italy & IT1 - IT10 & Latitude Longitude & Arbitrary aggregation & Municipality & 2001 & ISTAT\\
				\hline
				6 & Italy & ITp0 - ITp4 & Latitude Longitude & Arbitrary aggregation & Provincia & 2001 & ISTAT\\
				\hline
				7 & USA & USA0 - USA14 & Latitude Longitude & Arbitrary aggregation & County & 2000 & **\\
				\hline
	  	\end{tabular}
\end{table*}

\begin{figure*}
	\begin{center}
		\centerline{\includegraphics[scale=0.6]{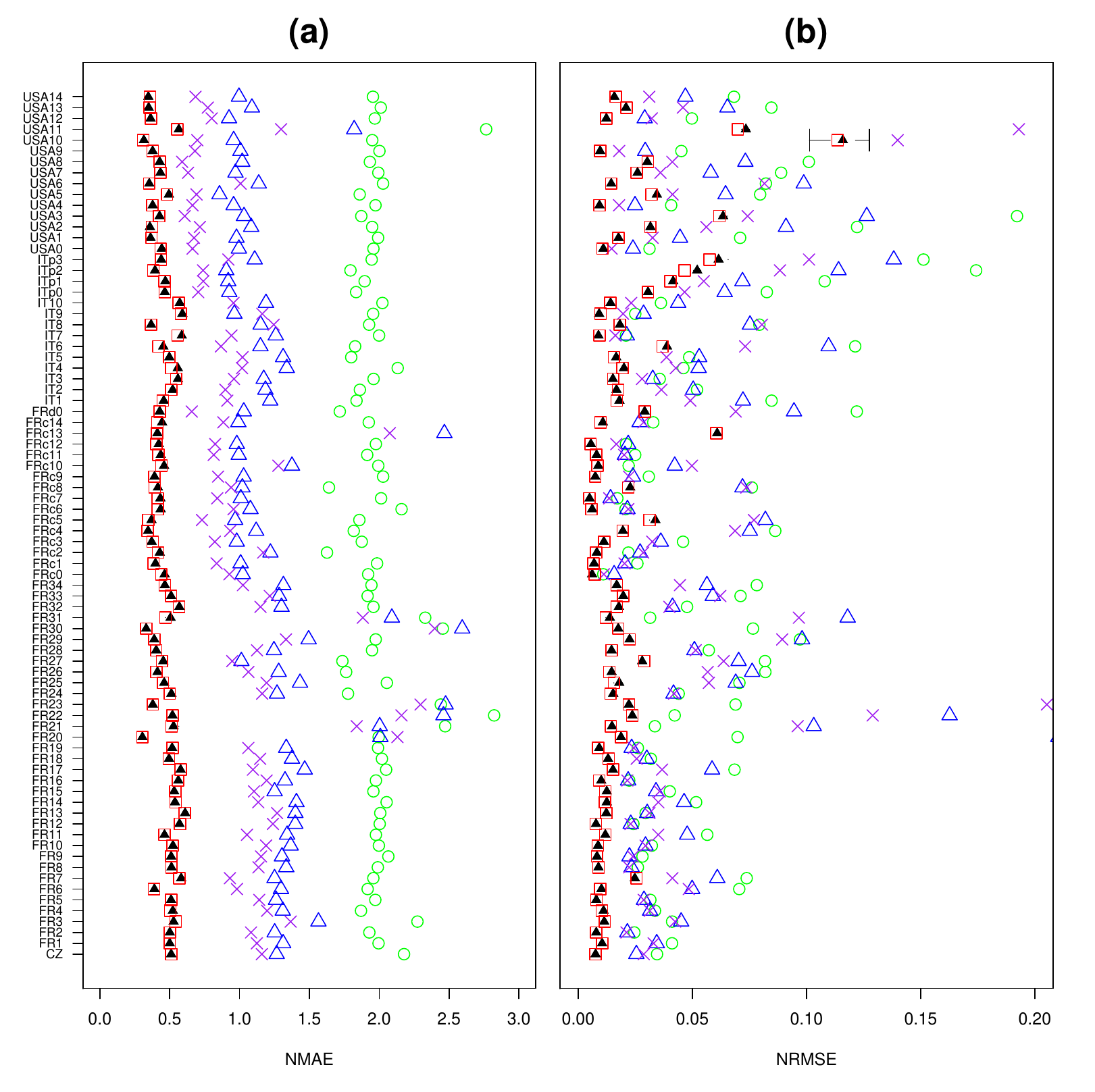}}
		\caption{Normalized Mean Absolute Error (a) and Normalized Root Mean Square 
		Error (b) for the 80 case-studies. 
		The red squares represent the errors obtained with the value of $\beta$ optimised from data on the case-study 
		network. Black plain triangles represent the average errors obtained with $\beta$ values estimated with the rule linking $\beta$ and the average 
		surface of the units obtain with the cross-validation; Dark bars represent the minimum and the maximum errors obtained with the estimated $\beta$ 
		but in most cases they are too close to the average to be seen. The green circles represent the errors obtained with the random model. The green 
		circles represent the errors obtained with the random model. The blue triangles represent the value obtained with the radiation model. The purple 
		cross represent the errors obtained with the modified version of the radiation model.
		\label{NMAENRMSENetwork}}
	\end{center}
\end{figure*}

\end{document}